\documentclass[12pt]{article}
\title{Riffles, ruffles, and the turning algebra}
\author{Peter G. Doyle and Dan Rockmore}
\date{PRELIMINARY Version dated 13 November 1995
\thanks{
Copyright (C) 1995 Peter G. Doyle.
Permission is granted to copy, distribute and/or modify this document
under the terms of the GNU Free Documentation License, 
as published by the Free Software Foundation;
with no Invariant Sections, no Front-Cover Texts, and no Back-Cover Texts.
}}
\newcommand{\rad}{\mbox{\tt rad}}
\newcommand{\gray}{\mbox{\tt gray}}
\newcommand{\even}{\mbox{\tt even}}
\newcommand{\post}{//}
\newcommand{\qed}{$\clubsuit$}
\newcommand{\rats}{{\bf Q}}
\newcommand{\cx}{{\bf C}}
\newcommand{\goesto}{\rightarrow}
\newcommand{\aut}{{\rm Aut}}
\newcommand{\semi}{\times}
\newcommand{\demi}{\times}
\newcommand{\cross}{\times}
\newcommand{\comp}{*}
\newcommand{\comptwist}{*}
\newcommand{\Radix}{\mbox{\tt Radix}}
\newcommand{\Gray}{\mbox{\tt Gray}}
\newcommand{\Riffle}{\mbox{\tt Riffle}}
\newcommand{\Ruffle}{\mbox{\tt Ruffle}}
\newcommand{\riffle}{{\mbox{\tt riffle}}}
\newcommand{\ruffle}{{\mbox{\tt ruffle}}}
\newcommand{\rising}{\mbox{\tt rising}}
\newcommand{\risingsequence}{\mbox{\tt risingsequence}}
\newcommand{\turning}{\mbox{\tt turning}}
\newcommand{\orientedruffle}{\mbox{\tt orientedruffle}}
\newcommand{\orientedrising}{\mbox{\tt orientedrising}}
\newcommand{\directedruffle}{\mbox{\tt directedruffle}}

\newcommand{\Sbar}{\bar{S}}
\newcommand{\direction}{ \{ \up,\down \} }
\newcommand{\up}{\mbox{\tt up}}
\newcommand{\down}{\mbox{\tt down}}

\newcommand{\hand}{\mbox{\tt hand}}
\newcommand{\cut}{\mbox{\tt cut}}
\newcommand{\colvec}[1]
  {{\left (
  \begin{array}{c}#1\end{array}
  \right )}}
\newcommand{\columnperm}[2]
  {
  {\left .
  \colvec{#1}
  \right |}
  \colvec{#2}
  }
\newcommand{\colperm}[2]
  {
  {\left (
  \begin{array}{c}#1\end{array}
  \right |}
  {\left .
  \begin{array}{c}#2\end{array}
  \right )}
  }
\newcommand{\permmatrix}{\mbox{\tt permmatrix}}
\begin{document}

\maketitle

\begin{abstract}
The \emph{rising algebra} is a subalgebra of the group algebra of the
symmetric group $S_n$, gotten by lumping together permutations having
the same number of rising sequences.
This well-known algebra arises naturally when studying riffle shuffles.
Here we introduce a number of other subalgebras that arise naturally
when stuffing `ruffles', which are like riffles except that 
after cutting the deck you turn over the bunch of cards that were
on the bottom.

This orphaned draft offers no context or motivation, and
uses idiosyncratic notation and terminology
that `seemed like a good idea at the time'.
We're making it available because it has
been cited in this form.
\end{abstract}

\section{To and fro}

\subsection{Natural order}

The theory of shuffling grows out of Jim Reeds's fundamental observation
that to understand the riffle shuffle, you have to look at it backwards.
Now, keeping straight the difference between $\sigma$ and $\sigma^{-1}$ is
a chore whenever you deal with permutations; having to try to keep everything
{\em backwards} is pretty near impossible (for us, at least).
To give ourselves a fighting chance, we have to write composition of functions
in the natural, left-to-right order.

{\bf WARNING.}
Throughout this paper we will compose functions in natural order:
\[
(\sigma \tau) [x] = \tau [ \sigma [x]]
.
\]

To try to minimize confusion,
we will use superscripts whenever possible,
so that
\[
x^{\sigma \tau} = (x^\sigma )^\tau
.
\]
We will also have occasion to use Wolfram's postfix notation, so that
\[
(x \post \sigma) \post \tau = x \post (\sigma \tau)
.
\]

\subsection{The permutation group}

Let $S_n$ denote the group of bijections from $\{1,\ldots,n\}$ to itself,
with functions composed in natural order:
\[
i^{\sigma \tau} = (i^{\sigma})^\tau
.
\]

As with any function defined on $\{1,\ldots,n\}$,
we can represent a permutation $\sigma$ as an $n$-tuple.
\[
(1^\sigma, \ldots, n^\sigma)
.
\]

We will adopt a variant of the cycle notation for permutations,
using $<i\ j>$ to denote the
transposition switching $i$ and $j$, and letting
\[
<i_1\ \ldots\ i_n> = <i_{n-1} i_n> <i_1\ \ldots\ i_{n-1}>
.
\]
(Don't forget: natural order!)

For example, if $n=3$, and $a,b$ are the standard `braid' generators
\[
a=<1\ 2> = (2,1,3)
,
\]
\[
b=<2\ 3> = (1,3,2)
,
\]
then
\[
ab = <1\ 2><2\ 3> = <1\ 3\ 2> = (3,1,2)
.
\]

This example demonstrates what appears to be 
a serious drawback of the $n$-tuple representation,
for while this representation `tells us where they went',
it doesn't {\em show} us.
To be more specific,
it seems very natural to represent the effects of $a$, $b$, and
$ab$ by drawing before-and-after diagrams:
\[
a = \frac{(1,2,3)}{(2,1,3)}
,
\]
\[
b = \frac{(1,2,3)}{(1,3,2)}
,
\]
\[
ab = \frac{(1,2,3)}{(2,3,1)}
.
\]
Here the numerator $(1,2,3)$ isn't conveying much information,
but if we omit it, the $n$-tuple that remains is the representation,
not of the permutation we're looking at, but of its inverse.
So it seems like the $n$-tuple representation
of a permutation is just backwards from what we want.  Where did we
go wrong?

Where we went wrong, of course, was in discarding the `numerator' of our
before-and-after diagram, which we should properly think of as a fraction.
If we interpret the fraction $\frac{\sigma}{\tau}$ to mean $\sigma \tau^{-1}$,
then everything is groovy:
\[
a = \frac{(1,2,3)}{(2,1,3)} = (2,1,3)
,
\]
\[
b = \frac{(1,2,3)}{(1,3,2)} = (1,3,2)
,
\]
\[
ab = \frac{(1,2,3)}{(2,3,1)} = (3,1,2)
.
\]

Moreover, if we rewrite $b$ as
\[
b = \frac{(1,2,3)}{(1,3,2)} = \frac{(2,1,3)}{(2,3,1)}
,
\]
then we get the very natural `braidlike' equations
\begin{eqnarray*}
a&=&
\frac{(1,2,3)}{(2,1,3)}
,\\
b&=&
\frac{(2,1,3)}{(2,3,1)}
,\\
ab&=&
\frac{(1,2,3)}{(2,3,1)}
.
\end{eqnarray*}

Of course we will have to be careful to remember that in the
formula
\[
\frac{\sigma}{\tau} = \sigma \tau^{-1}
,
\]
the $\tau^{-1}$ comes after (i.e., to the {\em right} of) the $\sigma$.
This is actually very natural,
if you consider that $\frac{\sigma}{\tau}$ is pronounced
`sigma divided by tau'.

This whole question of `Which came first, the sigma or the tau?' disappears if
we represent our $n$-tuples as column vectors, and transpose the fractions
$\frac{\sigma}{\tau}$ to $\sigma|\tau$.  Doing this yields the very congenial
equation
\[
ab =  \columnperm{1\\2\\3}{2\\1\\3}
\columnperm{2\\1\\3}{2\\3\\1}
= \columnperm{1\\2\\3}{2\\3\\1}
,
\]
which we abbreviate to
\[
ab =  \colperm{1\\2\\3}{2\\1\\3}
\colperm{2\\1\\3}{2\\3\\1}
= \colperm{1\\2\\3}{2\\3\\1}
.
\]
This vertical representation of permutations is particularly
appropriate in discussing shuffling, since it makes it easy to
visualize the cards as they appear in the deck.

The conventions and notations that we have adopted fit in well with
the representation of permutations as matrices.
Here we take our cue from the theory of Markov chains,
where a probability distribution is most conveniently represented as
a row vector $(p_1,\ldots,p_n)$ of positive numbers summing to 1.
A permutation $\sigma$ corresponds naturally to the Markov transition
matrix $\permmatrix[\sigma]$, where
\[
\permmatrix[\sigma]_{ij} = \delta_{\sigma[i],j}
.
\]
Multiplying our row vector $(p_1,\ldots,p_n)$ by $\permmatrix[\sigma]$
(on the right!) yields
\[
(p_1,\ldots,p_n)\, \permmatrix[\sigma]
= (p_{\sigma^{-1}[1]},\ldots,p_{\sigma^{-1}[n]})
,
\]
which fortunately turns out to be the effect of taking the quantities
$p_i$ and moving them from their initial position $i$ to position $\sigma[i]$.
Holding our breath, we check, and to our delight we find that,
with $a$ and $b$ as above,
\begin{eqnarray*}
&&
\permmatrix[a]\, \permmatrix[b]
\\&=&
\left (
\begin{array}{ccc}
0&1&0\\
1&0&0\\
0&0&1
\end{array}
\right )
\left (
\begin{array}{ccc}
1&0&0\\
0&0&1\\
0&1&0
\end{array}
\right )
\\&=&
\left (
\begin{array}{ccc}
0&0&1\\
1&0&0\\
0&1&0
\end{array}
\right )
\\&=&
\permmatrix[ab]
\end{eqnarray*}
All's well with the world!

\section{Actions and reactions}

Let $G$ be a group and $M$ a monoid.
(Actually, this whole discussion might go through when $G$ is only a monoid,
but we prefer to assume $G$
is a group---if only for alphabetical reasons---until
there is some good reason for generalizing to a monoid.)

Let $G$ act on $M$ on the right by automorphisms, so that
\[
(m^g)^h = m^{gh}
\]
and
\[
(mn)^g = m^g n^g
.
\]
When we need to refer to this action by name, we will attach this name to
the associated homomorphism
$\rho: G \goesto \aut[M]$.
Here $\aut[M]$ denotes the group of automorphisms with
natural-order composition, and $\rho$ is our default name for group actions.

Given an action $\rho$ of $G$ on $M$,
the {\em semidirect product} $G \semi_\rho M$ is the monoid
consisting of the set $G \cross M$ together with the composition law
\[
(g,m) (h,n) =
(gh, m^h n)
.
\]
To check the associative law, we note that
\[
((g,m)(h,n))(i,o)
=
(g,m)((h,n)(i,o))
=
(ghi,m^{hi} n^i o)
.
\]
$G$ is isomorphic to the submonoid $G \cross \{1\}$ of $G \semi_\rho M$.
The map $(g,m) \mapsto (g,1)$ is a monoid-homomorphism onto this submonoid;
its kernel is the `normal' submonoid $\{1\} \cross M$, which is isomorphic
to $M$.

Given an action $\rho$ of $G$ on $M$,
a map
$\gamma: M \goesto G$ is called a {\em reaction to $\rho$}
if
\[
\gamma[m] \gamma[n]
=
\gamma[m^{\gamma[n]} n]
,
\]
i.e., if
\[
(\gamma[m], m)
(\gamma[n], n)
=
(\gamma[m] \gamma[n], m^{\gamma[n]} n)
=
(\gamma[m^{\gamma[n]} n], m^{\gamma[n]} n)
.
\]
This means that the set
\[
\{(g,m): g = \gamma[m]\}
\]
is a submonoid of $G \semi M$.
As a set, the elements of this submonoid correspond naturally to the
elements of $M$, only the product in $M$ has been twisted through the
interaction with $G$.  We denote this new product by $\comptwist_\gamma$
(leaving $\rho$ to be inferred from context), so that
\[
m \comptwist_\gamma n = m^{\gamma[n]} n
,
\]
and we denote the monoid $M$ with product $\comptwist_\gamma$ by
\[
G \demi_\gamma M
.
\]
If we ever to have to call this something, we will call it the
{\em demisemidirect product} of $G$ and $M$ with respect to $\rho$ and
$\gamma$.

\section{Riffles and ruffles}

\subsection{The radix monoid}
Let $n$ be a positive integer, e.g.\ 52.
Denote by $\Radix_n$ the monoid with elements
$(a,(x_1,\ldots,x_n)): a \geq 1, 0 \leq x_i <n$,
and multiplication
\[
(a,(x_1,\ldots,x_n))(b,(y_1,\ldots,y_n))
=
(ab,(b x_1+y_1, \ldots, b x_n + y_n)
,
\]
or in simplified notation,
\[
\colvec{x_1\\ \vdots \\ x_n}_a
\colvec{y_1\\ \vdots \\ y_n}_b
=
\colvec{b x_1 +y_1\\ \vdots \\ b x_n + y_n}_{ab}
.
\]

We think of the elements of $\Radix_n$ as lists of digits
in the specified radix $a$;
we combine two lists entry-by-entry (in the natural order!),
interpreting the product of $\rad(x,a)$ and $\rad(y,b)$ as $\rad((x,y),(a,b))$,
a two-digit number in the hybrid radix $(a,b)$, where the first
digit $x$ is in radix $a$ and the second digit $y$ in radix $b$.

Note that if we represent $\rad(x,a)$ as the linear polynomial
$a X + x$, then the mixed-radix product of $\rad(x,a)$ and $\rad(y,b)$
corresponds to the composition (in natural order) of the corresponding
linear functions:
\[
(X \mapsto a X + x) (X \mapsto b X + y) = (X \mapsto ab X + bx + y)
.
\]

\subsection{The riffle monoid}
Now let $S_n$ act on $\Radix_n$ by permuting the list entries, and let
$\Radix_n$ react via the function $\riffle$ by interpreting the
entries in the list of digits base $a$
as portraying the effect of an $a$-handed riffle:
\[
\colvec{1\\1\\0\\1\\0}_2 \post \riffle
=
\colperm{1\\2\\3\\4\\5}{3\\4\\1\\5\\2}
,
\]
\[
\colvec{2\\2\\1\\0\\1}_3
\comptwist_\riffle
\colvec{1\\1\\0\\1\\0}_2
=
\colvec{1\\0\\2\\1\\2}_3
\comp_{\mbox{rad}}
\colvec{1\\1\\0\\1\\0}_2
=
\colvec{3\\1\\4\\3\\4}_6
,
\]
\begin{eqnarray*}
&&
\colvec{2\\2\\1\\0\\1}_3 \post \riffle
\colvec{1\\1\\0\\1\\0}_2 \post \riffle
\\&=&
\colperm{1\\2\\3\\4\\5}{4\\5\\2\\1\\3}
\colperm{1\\2\\3\\4\\5}{3\\4\\1\\5\\2}
\\&=&
\colperm{1\\2\\3\\4\\5}{4\\5\\2\\1\\3}
\colperm{4\\5\\2\\1\\3}{2\\1\\4\\3\\5}
\\&=&
\colperm{1\\2\\3\\4\\5}{2\\1\\4\\3\\5}
\\&=&
\colvec{3\\1\\4\\3\\4}_6 \post \riffle
.
\end{eqnarray*}

We call the monoid arising from this reaction
the {\em riffle monoid}:
\[
\Riffle_n = S_n \demi_{\riffle} \Radix_n
.
\]

\subsection{The Gray monoid}

As a variation on $\Radix_n$, we introduce $\Gray_n$, which is
to $\Radix_n$ as the Gray code is to binary.
Specifically, $\Gray_n$ has the same elements as $\Radix_n$, but
the new multiplication
\[
\colvec{x_1\\ \vdots \\ x_n}_a
\comp_{\mbox{gray}}
\colvec{y_1\\ \vdots \\ y_n}_b
=
\colvec{b x_1 + (x_1\ \even?\; y_1:b-1-y_1\\
\vdots \\ b x_n + (x_n\ \even?\; y_n:b-1-y_n}_{ab}
.
\]
Here we combine the Gray digits $\gray(x,a)$
and $\gray(y,b)$ by treating $(x,y)$ as a two-digit number
in the hybrid Gray base $(a,b)$, where the lower order Gray digit
runs alternately up and down, so that for example counting in Gray base
$(3,2)$ goes
\[
(0,0),(0,1),(1,1),(1,0),(2,0),(2,1)
.
\]

\subsection{The ruffle monoid}
To describe up-down riffles, or {\em ruffles},
we use the monoid $\Ruffle_n$, which we get by letting $S_n$ act as usual
on $\Gray_n$, and letting $\Gray_n$ react via the function $\ruffle$
by interpreting the
entries in the list of digits
as portraying the effect of an $a$-handed ruffle:
\[
\colvec{1\\1\\0\\1\\0}_2 \post \ruffle
=
\colperm{1\\2\\3\\4\\5}{5\\4\\1\\3\\2}
,
\]
\[
\colvec{1\\1\\2\\0\\1}_3
\comptwist_\ruffle
\colvec{1\\1\\0\\1\\0}_2
=
\colvec{1\\0\\1\\2\\1}_3
\comp_{\mbox{gray}}
\colvec{1\\1\\0\\1\\0}_2
=
\colvec{2\\1\\3\\5\\3}_6
,
\]
\begin{eqnarray*}
&&
\colvec{1\\1\\2\\0\\1}_3 \post \ruffle
\colvec{1\\1\\0\\1\\0}_2 \post \ruffle
\\&=&
\colperm{1\\2\\3\\4\\5}{4\\3\\5\\1\\2}
\colperm{1\\2\\3\\4\\5}{5\\4\\1\\3\\2}
\\&=&
\colperm{1\\2\\3\\4\\5}{4\\3\\5\\1\\2}
\colperm{4\\3\\5\\1\\2}{2\\1\\4\\5\\3}
\\&=&
\colperm{1\\2\\3\\4\\5}{2\\1\\4\\5\\3}
\\&=&
\colvec{2\\1\\3\\5\\3}_6 \post \ruffle
.
\end{eqnarray*}

This reaction yields the {\em ruffle monoid}:
\[
\Ruffle_n = S_n \demi_{\ruffle} \Gray_n
.
\]

\section{New algebras from old}
\subsection{Lumped monoids}
A function $\mu:M \goesto S$ from the monoid $M$ to an arbitrary set $S$
determines the equivalence relation $\equiv_\mu$, where $a \equiv_\mu b$
if and only if $\mu[a] = \mu[b]$.  We say that the function $\mu$ is
a {\em lumping}
if (the characteristic functions of) the
$\mu$-equivalence classes constitute a basis for a subalgebra of the monoid
algebra $\rats[M]$ (or $\cx[M]$, if you prefer).
(See Pitman
[?].)
Combinatorially, this amounts to requiring  that the $\mu$-equivalence
classes $[a]$ all be finite,
and that there exist {\em structure constants}
$C_{[a],[b],[c]}$ such that for any $a,b,c \in M$
there are exactly $C_{[a],[b],[c]}$ ways of writing $xy=c$ with
$x \in [a]$, $y \in [b]$.

\subsection{Do the right thing}
Let $M$ and $N$ be monoids, $\mu$  a lumping of $M$,
and $\nu$ a function on $N$ that we hope to show is a lumping.
We say that a homomorphism
$f:M \goesto N$ {\em does the right thing} if in the monoid algebra $\rats[N]$
the elements $\sum_{x \in [a]_\mu} f(x)$ belong to and span the subspace
spanned by (the characteristic functions of) the $\nu$-equivalence classes.
Conbinatorially, this amounts to requiring that there exist a matrix 
$D = \{D_{[a]_\mu,[b]_\nu}\}$ of what we might
call {\em restructure constants},
such that for any $a \in M$, $b \in N$
there are exactly $D_{[a]_\mu,[b]_\nu}$ ways of writing $f(x)=b$
with $x \in [a]_\mu$; 
in addition, the row-space of the matrix $D$ must contain
the standard basis vectors,
a requirement that in our examples will follow from the fact that
we can order the rows and columns of the matrix $D$
so that it becomes lower-triangular,
with non-zero entries on the diagonal.

{\bf Theorem.}
If $f:M \goesto N$ does the right thing with respect to
a lumping $\mu$ of $M$
and a function $\nu$ on $N$ then $\nu$ is a lumping of $N$.
\qed

\section{Shuffling and its algebras}

\subsection{Hand-equivalence and cut-equivalence}

In the monoids $\Riffle_n$, and $\Ruffle_n$ we can
lump elements together according to the value of the radix $a$.
Let's call the resulting equivalence relation $\hand$-equivalence,
since we are lumping together shuffles involving the same number of hands.
Note that the subalgebra yielded by the lumping $\hand$ is commutative:
Indeed, it is isomorphic to the monoid algebra of the natural numbers,
because every $ab$ riffle arises in one and only one way as an
$a$-riffle followed by a $b$-riffle (or vice versa).

Alternatively, we can refuse to identify two lists
unless in addition to sharing the same radix $a$, each base $a$ digit
occurs the same number of times in the second list as it does in the first.
This more discerning equivalence relation we call $\cut$-equivalence,
since now we are lumping together shuffles only if the cards
are cut and distributed among the $a$ hands in the same way.

\subsection{Rising sequences}

Given a permutation $\sigma \in S_n$, we cut the sequence $1,\ldots, n$
into subsequences called the {\em rising sequences} of $\sigma$ by dividing it
between $i$ and $i+1$whenever
$\sigma[i+1] < \sigma[i]$.
The number of rising sequences in $\sigma$
tells the minimum number of hands you need
in order to produce $\sigma$ as the result of a single riffle,
and the specific division into rising sequences tells where you have
to make the cuts in order to accomplish this.

The notion of rising sequences suggests
two equivalence relations on $S_n$.
We say that two permutations are
$\rising$-equivalent if they have
the same number of rising sequences,
and $\risingsequence$-equivalent
if in addition the rising sequences of the 
two permutations are exactly the same.

The map $\riffle$
does the right thing with respect to $\hand$ on $\Riffle_n$
and $\rising$ on $S_n$.
To verify this, we must check that the number of ways of
realizing a given permutation $\sigma$ as the result of an $a$-handed
shuffle depends only on the rising number of $\sigma$.
This fundamental observation about riffles is due to Bayer and Diaconis
[?]; the proof is a standard `stars-and-bars' argument.

Since $\riffle$ does the right thing,
$\rising$ is a lumping, and yields a commutative
subalgebra of the group algebra of $S_n$, which we call the
{\em rising algebra}.
(See
Bayer and Diaconis
[?],
Pitman
[?].)

The map $\riffle$ also does the right thing with respect to
$\cut$ and $\risingsequence$,
thus yielding the larger {\em rising sequence algebra}.
This second algebra is commonly called the `descent algebra'
(see
Bayer and Diaconis
[?],
Hanlon
[?];
we prefer to call it the `rising sequence algebra' because we feel this
is more in line with the general policy:
\begin{verse}
As you go through life make this your goal,\\
Watch the doughnut, not the hole!
\end{verse}

\subsection{Turning points and oriented rising sequences}

We say that a permutation has a turning point at $i, 1 \leq 1 <n$, 
if the graph of the permutation,
extended to map 0 to 0 and interpolated linearly to give a piecewise linear
mapping from $[0,n]$ to itself, has a local maximum or minimum at $i$.
Note that $1$ and $n$ are treated differently in this definition,
in that we call 1 a turning point if 2 ends up coming before 1, but we never
call $n$ a turning point.
The more symmetrical notion of `reduced' turning number
will be discussed later.

We say that two permutations are $\turning$-equivalent if they
have the same number of turning points.  (The stronger notion requiring
that in addition they have exactly the same set of turning points
coincides with the $\risingsequence$ relation, 
so we can ignore it.)

While the rising number of the identity permutation is 1,
the turning number of the identity
is 0, and in general the turning number of a permutation is 1 less that we
would hope and expect.
This anomaly stems from the fact that
while rising sequences begin and end at the interstices between
consecutive integers in the domain of the permutation,
turning points occur at the integers themselves.
To get a more felicitous analog of rising sequences,
we must look not at permutations, but at {\em oriented permutations}
(also called `signed permutations').
In card-shuffling terms, an oriented permutation keeps track
of the way the cards are facing as well as their order.
We divide the sequence $1,\ldots,n$ into subsequences,
called {\em oriented rising sequences},
according to the cuts we would need to make in order to achieve
the specified arrangement in a single ruffle with the minimum number of hands.
Note that some of these subsequences may have length 0, though no two
of them in a row will have length 0.
The maximum oriented rising number
is $2n$:
To turn over each card in place with a since ruffle, you need $2n$
hands.

The turning number of a permutation can then be viewed as one less than
the minimum oriented rising number of an oriented permutation that
reduces to the given permutation when the orientations of the cards
are ignored.

\subsection{The turning algebra and the oriented rising algebra}
The map $\ruffle: \Ruffle_n \goesto S_n$ does the right thing
with respect to $\hand$ and $\turning$.
As you would expect,
the best way
to see that $\ruffle$ does indeed do the right thing
is to factor it through
the group $\Sbar_n$ of oriented permutations.
Writing
\[
\ruffle[m] = \orientedruffle[m] \post \pi_{\Sbar_n,S_n}
,
\]
we observe first that $\orientedruffle$ does the right thing with respect
to $\hand$ on $\Ruffle_n$
and $\orientedrising$ on $\Sbar_n$.
Thus $\orientedrising$ is a lumping.
Then we observe
that $\pi$ does the right thing with respect to $\orientedrising$
and $\turning$, so $\turning$ is a lumping.
Verifying that these two maps do indeed do the right thing
involves showing that the number of ways of obtaining a given
oriented permutation as the result of an $a$-ruffle depends only
on the oriented rising number, and the number of ways of obtaining
a given permutation as the result of an oriented permutation with specified
oriented rising number depends only on the turning number.
As in the case of riffles, these verifications involve elementary
counting arguments.

Thus we obtain a commutative subalgebra of the group algebra
of $S_n$, the {\em turning algebra},
by way of a commutative subalgebra of the group algebra
of $\Sbar_n$, the {\em oriented rising algebra}.

\subsection{The reduced turning number}
The {\em  reduced turning number} of a permutation differs from
the turning number in that it refuses to recognize a turning point at 1.
Thus the identity permutation and the permutation $\Delta$ that reverses
$1,\ldots,n$ (turning over the deck) both have reduced turning number 0.

To treat the reduced turning number with the machinery we have developed,
we need a version of ruffling where instead of always turning over the
odd-numbered piles, we turn over the odds or the evens depending on a
specified direction.
Thus we replace $\Gray_n$ with
$\direction \cross \Gray_n$.
We extend the action of $S_n$ by having it leave
the direction alone, and we let
$\direction \cross \Gray_n$.
react by interpreting the list of digits as the effect of
a regular (`up-down') ruffle or a reverse (`down-up') ruffle
according to the value of the direction:
\[
(\up,(2,(1,1,0,1,0))) \post \directedruffle
=
\frac{(1,2,3,4,5)}{(5,4,1,3,2)}
,
\]
\[
(\down,(2,(1,1,0,1,0))) \post \directedruffle
=
\frac{(1,2,3,4,5)}{(3,4,2,5,1)}
.
\]
This reaction yields the {\em directed ruffle monoid}:
\[
S_n \demi_{\directedruffle} (\direction \cross \Gray_n)
\]

The map $\directedruffle$
does the right thing (here again by way of $\Sbar_n$), and we get the
a commutative subalgebra
of the group algebra of $S_n$,
the {\em reduced turning algebra},
by way of the corresponding
commutative subalgebra of the group algebra of $\Sbar_n$.

\section{Bijective correspondences}

The existence of the rising algebra is equivalent to the fact
that any two permutations with $k$ rising sequences arise in the same
number of ways as products of permutations with $i$ and $j$ rising sequences.
Moreover, in this case there is a natural bijection between the sets
of factorings of the two permutation (and also between these factorings
and the factorings where the roles of $i$ and $j$ are reversed).
The existence of these bijections follows on general principles from
the fact that the matrix of restructure constants associated with
the map $\riffle$, which induces the rising algebra from the
combinatorially trivial $\hand$-equivalence subalgebra of the riffle
algebra, can be written in lower-triangular form with
{\em 1's on the diagonal}.
The same goes for the oriented rising algebra.
However, in the case of
the turning algebra,
the diagonal restructure constants are powers of two.
Of course, since the sets of factorings in question are the same size,
there will of still exist bijections.
However, there will no longer be any reason to expect that there will
exist any natural bijections
(whatever that might mean).
Thus it appears that from a combinatorial point of view,
the turning algebra is essentially more
complicated
than
the oriented rising algebra (from which it arises by lumping)
and the rising algebra.

\end{document}